%Paper: funct-an/9212003
%From: Alan Hopenwasser <AHOPENWA@UA1VM.UA.EDU>
%Date: Fri, 18 Dec 92 11:47:31 CST

%AMS-TEX Version 2.1
\input amstex
\documentstyle{amsppt}
\magnification =1200
\pagewidth {6.5 true in}
\pageheight {8.5 true in }
\hcorrection {-.2 true in}
\vcorrection {+.1 true in}
\baselineskip = 14pt
\parindent 20 pt
\parskip = 3pt plus .5pt minus .5pt
\define\ds{\displaystyle}

\define\lb{\lbrack}
\define\rb{\rbrack}
\define\braks#1{\lb \,#1\,\rb}
\define\nospacebraks#1{\lb #1\rb}
\redefine\H{\Cal H}
\define\K{\Cal K}
\define\N{\Cal N}
\define\kn{\K_{\N}}
\redefine\D{\Cal D}
\define\bh{\Cal B (\Cal H)}
\define\cbb{\Bbb C}

\define\nbb{\Bbb N}

\define\system#1#2#3#4{T_{#1} @> \phi_1 >> T_{#2}
@> \phi_2 >> T_{#3} @> \phi_3 >>
\dots @>>> #4}

\define\Asystem#1#2#3#4{A_{#1} @> \phi_1 >> A_{#2}
@> \phi_2 >> A_{#3} @> \phi_3 >>
\dots @>>> #4}

\define\alg{\Cal A\text{\sl lg\/}\,}

\define\algn{\alg (\Cal N)}

\define\plstar{${}^*$}
\define\cstar{$\text{C}^*$}
\topmatter
\title
Compression Limit Algebras
\endtitle
\author
Alan Hopenwasser \\
Cecelia Laurie
\endauthor
\affil
University of Alabama
\endaffil
\address
Department of Mathematics, University of Alabama,
Tuscaloosa, AL 35487
\endaddress
\email
ahopenwa\@ua1vm.ua.edu and claurie\@ua1vm.ua.edu
\endemail
\thanks
The authors would like to thank Vern Paulsen for
helpful suggestions on the content of this paper.
\endthanks
\endtopmatter
\document
\head
I. Introduction
\endhead

While there is now a large and rapidly growing literature on
the study of direct limits of subalgebras of
 finite dimensional \cstar-algebras,
the focus in almost every paper on the subject has been on systems
with embeddings which have \plstar-extensions to the
generated \cstar-algebra.
One notable exception is the paper by Power \cite{P1}.
The purpose of this note is to extend the study of
systems which are not \plstar-extendible and to
produce examples which illustrate some new phenomena.
\par
When all the embeddings in a  direct system
are \plstar-extendible, then the limit algebra
is, in a natural way, a subalgebra of an AF
\cstar-algebra. On the other hand,
 if the embeddings are not
\plstar-\nolinebreak extendible, then it is no longer
a~priori obvious that the limit algebra is
an operator algebra.  (Initially, the direct
limit must be taken within the category of
Banach algebras.)
Even if the limit algebra is an operator algebra,
its image under some representations may generate
a \cstar-algebra which is not approximately
finite (as happens in the situation studied
by Power).  For systems with the type of
embeddings which we consider in this paper,
the limit algebra will always be an operator
algebra, a fact which can easily be seen
with the aid of the abstract characterization
of operator algebras by Blecher,
 Ruan, and Sinclair  \cite{BRS}.
We also produce a representation of the limit
algebra which is ``natural'' in the sense
that the \cstar-envelope  of the image algebra
 (as defined by
Hamana \cite {H}) is isomorphic to the
\cstar-algebra generated by the image algebra.

\head
II. Compression Embeddings
\endhead

The key observation behind the choice of embeddings
under investigation is that if $A$ is a CSL-algebra
and if $p$ is an interval from the lattice of
invariant projections for $A$, then the mapping
$x \mapsto pxp$ is an algebra homomorphism.
Since we are interested in direct limits
of finite dimensional algebras, we shall assume
that every CSL-algebra is finite dimensional.
This implies that, with respect to a suitable
choice of matrix units, $A$ is a subalgebra
of some full matrix algebra $M_n$ which contains
the algebra  $D_n$  of diagonal matrices.  Such
algebras go under a variety of names:
 {\it poset algebras\/},
{\it incidence algebras\/}, or
{\it digraph algebras\/}; henceforth,
we       use the term digraph algebra.
\par
The invariant projections of a digraph algebra are
all projections in the diagonal  $D_n$; in particular,
they all commute with one another.  If $e$ and $f$
are two invariant projections such that $e \le f$,
then $p = f-e$ is called an {\it interval\/}
 from the lattice.
(The term {\it semi-invariant\/} projection
 is also sometimes used.)
\par
Since we want our embeddings to be unital, a
{\it compression\/}
        will be the mapping which takes $x$ to the
restriction of $pxp$ to the range of $p$,
where $p$ is an interval from the lattice.  The image
algebra under the compression will be viewed as a
subalgebra of a full matrix algebra, usually of rank
smaller than the rank of the original containing
matrix algebra.
\definition{Definition}  Let $A \subseteq M_n$ be a
digraph algebra.
A {\it compression embedding\/}           is a direct sum
of compression mappings on $A$ subject to the
proviso that at least one summand is the
identity mapping.
\enddefinition
The purpose of the assumption that at least one
summand is the identity mapping is to ensure
that the embedding  is completely isometric.
(Compressions are completely contractive, but not
completely isometric except in the trivial case
of the identity mapping.)
\definition{Definition}  A {\it compression limit
algebra\/}       is the limit of a direct system
$$
\Asystem 123A
$$
in which each embedding            is a unital
compression embedding.
\enddefinition
\par
A compression embedding of $A_k$ into $A_{k+1}$
is a unital map of $A_k$ onto its range; we are
assuming then that the unit of $A_{k+1}$ is equal to
the unit of
the subalgebra  $\phi_k(A_k)$.  (This avoids
unwanted degeneracies.)  The special case in
which each summand is compression to the identity
is just the familiar case of a standard embedding.
In this special case, compression embeddings are
\plstar-extendible; in general they will not be.
\par
In this note, we shall focus on a special class
of compression algebras: direct limits of full
upper triangular matrix algebras,
$$
\system {n_1}{n_2}{n_3}A.
$$
The non-\plstar-extendible embeddings studied
by Power \cite{P1} require non-trivial
cohomology for the digraph algebras in the system;
hence nest algebras can never appear as finite
dimensional approximants.  The systems studied
here have one other new feature: the only
constraint on $n_k$ and $n_{k+1}$ is that
$n_k < n_{k+1}$; no divisibility is required.
\par
For each $k$, one of the summands of $\phi_k$
is the identity compression; let $q_k$ denote the
support interval for this summand from the nest of
invariant projections for $T_{n_{k+1}}$.
Also, let $\psi_k$ denote the compression of
$T_{n_{k+1}}$ onto $q_k$.  We may, in  a natural way,
view $\psi_k$ as an algebra homomorphism from
$T_{n_{k+1}}$ onto $T_{n_k}$.  It is clear that
$\psi_k \circ \phi_k = id$, for all $k$.  Thus
compression embedding systems are structured in
the sense introduced by Larson \cite{L};
in particular, compression limit algebras are all
structured Banach algebras.
\par
Compression embeddings satisfy one additional
property: they map matrix units to sums of
matrix units.  This implies that they are
regular embeddings, in the sense used
by Power \cite {P2, section 4.9}.
\par
\remark{Remark}
Compression embeddings can be placed into a broader
context.  If $C = \braks {c_{ij}}$
 is an $n \times n$ matrix, the
mapping $\phi_C \: M_n \to M_n$ given by
$\phi_C \braks {a_{ij}} = \braks  {c_{ij}a_{ij}}$
is known as a {\it Schur\/} mapping.
It is easy to check that $\phi_C$ is an algebra
homomorphism if, and only if, $C$ satisfies the
cocycle condition: $c_{ik} = c_{ij}c_{jk}$.
Furthermore, $\phi_C$ is unital if, and only if,
all $c_{ii} = 1$.
By taking a direct sum of algebra homomorphisms
each of which is a compression composed with a
Schur mapping, we can obtain embeddings more
general than compression embeddings.  (Compression
embeddings arise by insisting that appropriate
$c_{ij} = 1$.)  We defer the study of these
more general systems to another time.
\endremark

\head
III. Representations
\endhead

Let
$$
\system {n_1}{n_2}{n_3}A
$$
be a direct system with compression embeddings.
Since each $\phi_k$ is a unital complete isometry,
we may identify each $T_{n_k}$ with a unital
subalgebra of $A$; the subalgebras so obtained are
nested and the closure of their union is $A$.
Each subalgebra carries a matricial norm and the
sequence of norms is  compatible with the nesting;
consequently, $A$ has a matricial norm.
It is easy to see that this matricial norm satisfies
the axioms of Blecher, Ruan and Sinclair
\cite {BRS}, so $A$ is an (abstract) operator
algebra, i.e., there is a Hilbert space $\H$
and a completely isometric unital algebra
homomorphism $\rho$  mapping $A$ into $\bh$.
This argument remains valid for systems of
digraph algebras and for systems with
the more general Schur embeddings,
provided that the embeddings
are complete isometries.
\par
Our primary interest, however, is to obtain
an explicit representation which is, in an
appropriate sense, canonical.  To that end,
let $\H$ denote a Hilbert space with basis
$\{e_n\}$.  This index set will depend
on the specific direct system, but will
always be either $\Bbb Z$ or $\{\,n\: n \le b\,\}$
for some integer $b \ge 0$ or
$\{\,n\:n \ge a\,\}$ for some integer
$a \le 0$.
For each $k$, let $p_k$ denote the projection
on the closed linear span of $\{\,e_n \: n \le k\,\}$
and let $\N$ denote the nest consisting
of the projections $p_k$ together with 0 and $I$.
$\algn$ denotes the nest algebra associated with $\N$.
The representation which we will construct will map
$A$ to a weakly dense subalgebra of $\algn$.
\par
For each $k$, we may write $\phi_k =
\psi^k_a \oplus \dots \oplus \psi^k_b$, where
each $\psi^k_j$ is a compression map.  At least one of
these compression maps must be the identity; in order to
keep track of a selection of identity compressions, arrange
the indexing so that $\psi^k_0 = id$ for all $k$.
Thus, we will always have that the first index $a$
is non-positive and the last index $b$ is non-negative.
\par
Now suppose that $\phi\: T_{n_k} \to T_{n_{k+1}}$ is a
compression embedding and that $\psi$ is a single
compression to an interval from the invariant
projection lattice of $T_{n_{k+1}}$.
It is easy to check that $\psi \circ \phi$ is
a sum of compressions to intervals from the lattice
of $T_{n_k}$.  Note that the intervals associated with
$\psi \circ \phi$ may include some which were not associated with
$\phi$ itself.
\par
A consequence of the observation above is that a
composition of two compression embeddings is
again a compression embedding.  Since we need to
consider multiple compositions, let $\phi_{jk}$
denote the composition $\phi_k \circ \dots \circ \phi_j$.
Then $\phi_{jk}$ can be written in the form
$\ds \sum_a^b{}^{\oplus}\eta_i$, where each $\eta_i$
is a simple compression.  Furthermore, since
$\psi_0^n = id$ for each $n$, we may arrange the
indexing so that $\eta_0 =id$.  Thus, $a \le 0$ and
$b \ge 0$.  The indices $a$ and $b$ depend, of course,
on $j$ and $k$, but that dependence is suppressed
in the notation.
\par
Observe also that the string of $\eta$'s
        in the summation for $\phi_{jk}$
appears as a substring in the expression
for $\phi_{j,k+1}$.  The reason, of course,
is that $\psi_0^{k+1} = id$.  Furthermore,
this substring includes the $\eta_0 = id$
term.  Thus, as we increase $k$, we change
the summation for $\phi_{jk}$ by adding terms at
one or both ends of the sum.  By increasing $k$
indefinitely, we obtain a mapping $\ds \rho_j =
\sum{}^{\oplus}\eta_i$.
The index set for this sum will either be
$\Bbb Z$ or a set of the form
$\{\,i\:i \ge a\,\}$ or of the form
$\{\,i\:i \le b\,\}$ for some $a \le 0$
or $b \ge 0$.
\par
An index set consisting of integers greater
than or equal to some $a$ will occur
precisely when all but finitely many of
the embeddings $\phi_k$ have $\phi_0^k = id$
as the  first term; the other singly infinite index
set arises when the distinguished identity summand
occurs as the last term in all but finitely many
embeddings.  In either case, we can change the
presentation by deleting finitely many embeddings
to arrange that the index set for the $\eta$'s
is either the non-negative integers or the
non-positive integers, as appropriate.
In either of these cases, we may now view
$\rho_j$ as a representation on the Hilbert
space $\H$ with respect to the appropriate basis.
The doubly infinite case still presents an
ambiguity: we need to ``anchor'' the representation
$\rho_j$ with respect to the basis $\{\,e_n \:n \in \Bbb Z\,\}$.
\par
This ``anchoring'' may be done as follows.  Let the
index set for the matrices in $T_{n_1}$ be
$\{\,0,1,\dots,b\,\}$, so that $f_{00}^1$ is
the matrix unit which has a 1 in the upper left
corner and zeros elsewhere.
Identify $\psi_0^1$ with the mapping
$0 \oplus \dots \oplus 0 \oplus \psi_0^1 \oplus 0
\oplus \dots \oplus 0$, a mapping of $T_{n_1}$
into $T_{n_2}$.  (In other words, replace all
terms in $\phi_1$  except for the distinguished
copy of the identity by a 0 mapping of appropriate
dimension.)
The image of $f^1_{00}$ under $\psi_0^1$ is now
a matrix unit in $T_{n_2}$; arrange the index set
for $T_{n_2}$ so that this matrix unit is
$f_{00}^2$.  By iterating this procedure, we
can index all the $T_{n_k}$ so that the
successive images of the matrix unit $f_{00}^1$
under the distingushed identity compressions
will be $f_{00}^k$.  In the doubly infinite case,
the indexing of the $T_{n_k}$'s will begin with
negative integers and end with positive integers
each of arbitrarily large magnitude for sufficiently
large $k$.
\par
It is now clear how to define the representation $\rho_k$:
the image of the distinguished matrix unit $f_{00}^k$
in the component $\eta_0$ should be the one dimensional
projection onto the span of the basis vector $e_0$.
For each $k$, $\rho_k$ is   a completely isometric representation
         of $T_{n_k}$ into $\algn$ acting on the Hilbert
space $\H$.   Furthermore,
 the following diagram commutes:
$$
\CD
T_{n_k} @>\phi_k>> T_{n_k+1} \\
@VV\rho_kV  @VV\rho_{k+1}V \\
\algn  @>id>> \algn
\endCD
$$
This system of representations induces a
completely isometric representation $\rho$
from the direct limit  $A$  into $\algn$.
\remark{Remark}
Observe that the image $\rho(A)$ is weakly (or
$\sigma$-weakly) dense in $\algn$.  To see this
it is enough to note that,
given $a \le 0$ and $b \ge 0$,
 the image contains
matrices with arbitrarily specified entries in
the locations $ij$ for $a \le i \le j \le b$.
To accomplish this,  choose $k$ sufficiently
large, select an appropriate matrix in $T_{n_k}$, and
inspect the image of that matrix under the representation $\rho$.
\endremark
\par
In the special case in which each $\phi_k$ is a direct
sum of identity mappings, i.e., $\phi_k$ is a
standard embedding, the representation which we have
constructed is just the representation of the standard
limit algebra introduced by  Roger Smith.  For a description
of this representation, see \cite{OP}.
Depending on the locations of the distinguished copy
of the identity, the representation will act on a
Hilbert space with a basis indexed by the integers,
by the non-negative integers, or by the non-positive
integers.  The choice here is arbitrary, but since
the representations in this special case are all
\plstar-extendible, the \cstar-algebras generated
by the images are all isomorphic.  (Indeed, they
are all the UHF algebra associated with the appropriate
supernatural number.)
\par
In the general case, it is possible to have two
completely isometric representations, $\rho$ and
$\sigma$, of $A$ with $C^*(\rho(A))$ not
isomorphic to $C^*(\sigma(A))$.  It is desirable, then,
 to find a representation for which the \cstar-algebra
generated by the image is the \cstar-envelope in
the sense of Hamana.  We shall show below that $\rho$,
as constructed above, has just this property.  It is
in this sense that $\rho$ is canonical.
\par
First, we provide a brief review of the idea of a
\cstar-envelope.  The appropriate setting for this
is actually the category of
 unital operator spaces with unital
complete  order injections, but we only
need to deal with operator algebras so we will
restrict the discussion to that domain.
\par
If $A$ is an operator algebra, then a \cstar-{\it extension}
of $A$ is a \cstar-algebra $B$ together with a unital
complete order injection $\rho$ of $A$ into $B$ such
that $C^*(\rho(A)) = B$.
A \cstar-extension $B$ is a \cstar-{\it envelope} of $A$
provided that, given any operator system, $C$, and any
unital completely positive map $\tau \: B \to C$,
$\tau$ is a complete order injection whenever
$\tau \circ \rho$ is.  Hamana \cite{H}
proves the existence and uniqueness (up to a suitable
notion of equivalence) of \cstar-envelopes; further, he
shows that the \cstar-envelope of $A$ is a minimal
\cstar-extension in the family of all
 \cstar-extensions of $A$.
\par
\define\Silov{\v Silov}
After proving the existence of \cstar-envelopes, Hamana
then uses this to prove the existence of a \Silov \ boundary
for $A$.  The \Silov \ boundary is a generalization to
operator spaces of the usual notion of \Silov \ boundary
from function spaces; it was first developed by Arveson
in \cite{A1}.  Here is  the appropriate definition.
Let $B$ be a \cstar-algebra and let $A$ be a unital
subalgebra such that $B = C^*(A)$.  An ideal  $J$
in $B$ is called a {\it boundary ideal} for $A$
if the canonical quotient map $B \to B/J$ is completely
isometric on $A$.  A boundary ideal which contains
every other boundary ideal is called the {\it \Silov \
boundary\/} for $A$.
\par
Hamana shows that if $B$ is a \cstar-extension
for $A$, then the \cstar-envelope of $A$ is
isomorphic to $B/J$, where $J$ is the \Silov \
boundary for $A$.  In particular, $B$ is the
\cstar-envelope for $A$ if, and only if, the
\Silov \ boundary is 0.  This last fact is the
one which we shall use to show that $C^*(\rho(A))$
is the \cstar-envelope of the image $\rho(A)$ for
the representation $\rho$ defined above.
\par
The first step needed to accomplish this goal
is to compute the \cstar-algebra generated by
              each $\rho_k(T_{n_k})$.  The
following simple lemma is helpful:
\proclaim{Lemma} Let $T_n$ be the $n \times n$ upper
triangular matrices acting on $\Bbb C^n$ and let
$p$ and $q$ be distinct intervals from the nest of
invariant projections for $T_n$.  Then
$$
C^*(\{\,pap \oplus qaq \: a \in T_n\,\}) =
\Cal B(p\Bbb C^n) \oplus \Cal B(q\Bbb C^n).
$$
\endproclaim
\demo{Proof}
Here, we view $pap$ and $qaq$ as being restricted to
the ranges $p\Bbb C^n$ and $q\Bbb C^n$ of $p$ and $q$
respectively.
Since $p$ and $q$ are distinct, there is an element
$a$ of $T_n$ such that
 one of $pap$ and $qaq$, say
$pap$, is non-zero while the other is 0.  So we have
an element of the form $b \oplus 0$ in
$C^*(\{\, pap \oplus qaq \: a \in T_n\,\})$.
{}From this and the fact that the \cstar-algebra
generated by $T_n$ is $M_n$, it follows that
$  \Cal B (p\Bbb C^n) \oplus 0    \subseteq
C^*(\{\, pap \oplus qaq \: a \in T_n\,\} $\unskip.
This in turn implies that
                 $ 0  \oplus \Cal B(q\Bbb C^n)
\subseteq
C^*(\{\, pap \oplus qaq \: a \in T_n\,\} $ and hence that
$  \Cal B (p\Bbb C^n) \oplus \Cal B(q\Bbb C^n)
\subseteq
C^*(\{\, pap \oplus qaq \: a \in T_n\,\} $.
The reverse containment is evident.
\enddemo
By using the obvious extension of this lemma to multiple
direct sums (including countable direct sums), we can
describe the \cstar-algebra, $C^*(\rho_k(T_{n_k}))$, both
as  a subalgebra of $\bh$ and as an abstract (finite-dimensional)
\cstar-algebra.  Indeed, write $    \rho_k = \sum^{\oplus}
\eta_i$ and, for each $i$, let $q_i$ be the interval in the
nest for $T_{n_k}$ to which $\eta_i$ is a compression.
Then
\define\crk{C^*(\rho_k(T_{n_k}))}
\define\crj{C^*(\rho_j(T_{n_j}))}
\define\crho{C^*(\rho(A))}
$$
\crk = \{\, (b_j) \in \sum{}^{\oplus}
      \Cal B (q_j \Bbb C^{n_k})\:  b_i = b_j
     \text{ whenever } q_i = q_j \,\}.
$$
If $r_1, \dots , r_m$ is the list of distinct intervals
to which the $\eta_i$ are compressions,  then
$$
\crk \cong {\sum\limits_{j=1}^m}{}^{\oplus}
   \Cal B (r_j \Bbb C^{n_k}).
$$
We shall need both the isomorphism class of $\crk$ and
its actual expression as an operator algebra
 acting on $\H$.
Also, since $\crho$ is the closure of
the union of the $\crk$, we have proven the
following:
\proclaim{Proposition} Let $A$ be the direct limit
of a system of full upper triangular matrix algebras
with compression embeddings and let $\rho$ be the
representation of $A$ defined above.  Then
$\crho$ is an AF \cstar-algebra.
\endproclaim
Our main result is the following theorem.
\proclaim{Theorem} Let $A$ be the direct limit
of a system of full upper triangular matrix algebras
with compression embeddings and let $\rho$ be the
representation of $A$ defined above.  Then
$\crho$ is the   \cstar-envelope of $\rho(A)$.
\endproclaim
\demo{Proof}
In order to prove that $\crho$ is the
\cstar-envelope of $\rho(A)$, it is sufficient
to show that the \Silov\ boundary of $\rho(A)$
is 0.  Since the \Silov\ boundary is the largest
boundary ideal, we need merely show that if
$J$ is a non-zero ideal then $J$ is not a boundary
ideal.  This requires proving that the canonical
quotient map $\crho \to \crho /J$ is
not completely isometric when restricted to $\rho(A)$.
We shall, in fact, prove that $\rho(A) \cap J \ne \emptyset$;
the quotient map is not even isometric.
\par
For convenience,
let $B_k = \crk$ and $B = \crho$.  Each $B_k$ is
isomorphic to a direct
sum of full matrix algebras, one for each compression
map which appears in the expression for~$\rho_k$.
Since the identity must appear, one summand must be
$M_{n_k}$.  This summand is the largest rank summand and
it appears one time only.  $M_{n_k -1}$ appears at most
two times, reflecting the fact that the nest for $T_{n_k}$
has exactly two interval projections of rank $n_k -1$.
More generally, $M_{n_k -j}$ appears at most $j+1$ times
in the expression for $B_k$.
\par
Thus, the isomorphism class of $B_k$ is
$M_{n_k} \oplus M_{t_1} \oplus \dots
\oplus M_{t_s}$, where $t_1, \dots t_s$
are integers less than $n_k$.
In the Bratteli diagram for $B$, the $k^{\text{th}}$
level has a node for each summand in the isomorphism
class of $B_k$.
Next, we need the following observation: given $k$ and
a summand $M_d$ in the isomorphism class of $B_k$,
there is an integer $j >k$ such that $M_d$ partially
embeds into the summand $M_{n_j}$ in the isomorphism
class of $B_j$.
\par
The observation is verified by inspecting the chain of
finite dimensional \cstar-algebras acting on
$\H$. Now, $\rho_k(T_{n_k})$ is an infinite sum whose
terms are selected from finitely many compressions
of $T_{n_k}$.  Consequently, there exist integers
$a \le 0$ and $b \ge 0$ such that if $p$ is the
projection onto the linear span of $\{\,e_i \:
a \le i \le b\,\}$, then $p$ is reducing for
$\crk$ and every compression, and in particular
the one corresponding to $M_d$, appears in the
restriction of $\crk$ to $p$.  Then, by the way
in which the representations are constructed,
there is an integer $j >k$ so that the $\eta_0 = id$
term in $\rho_j$ acts on a subspace of $\H$ which
includes the range of $p$.  (This was the reason
for ``anchoring'' the $\rho_j$ so that the $\eta_0$
terms act on the vector $e_0$ and that the support space
for the $\eta_0$ terms ``grows away'' from $e_0$.)
This shows that in the abstract Bratteli diagram,
each node eventually partially
 embeds into a node corresponding
to the identity summand (the $M_{n_j}$ node).
\par
Now let $J$ be a non-zero closed two-sided ideal
in $B$.  By a result in \cite{B}, $J$ is the
closure of the union of the $J \cap B_k$.  In particular,
there is a positive integer $k$ such that $J \cap B_k \ne 0$.
Since $B_k$ is isomorphic to a finite direct sum
of full matrix algebras, $J \cap B_k$ is isomorphic to a
direct sum of some of those algebras, with 0's
as the remaining summands.  Let $M_d$ be one
of the non-zero summands appearing in $J \cap B_k$.  Let
$j >k$ be an integer such that the $M_d$ term
partially embeds into $M_{n_j}$.  It now follows
that $M_{n_j}$ is one of the non-zero summands for
$J \cap B_j$.
By utilizing the support projection for the $M_{n_j}$
term (a reducing projection for $\crj$), we see that
the subalgebra of $\crj$ isomorphic to $M_{n_j}$
is contained in $J$.
This subalgebra consists of all sequences $(b_i)$
in $\crj$ for which $b_i =0$ whenever the corresponding
interval $q_i$ is not the identity and the remaining
$b_i$ are all equal.
\par
The proof is now completed by observing that there
is an element of this subalgebra which lies in
$\rho_j(T_{n_j})$.  Indeed, let $v$
 be the matrix unit in $T_{n_j}$ which
has a 1 in the extreme upper right corner
and zeros elsewhere.   Then
the compression of $v$ to any interval $q$ other than
the identity is zero.  Thus, $\rho_j(v) \in J \cap \rho(A)$,
and the proof is complete.
\enddemo
\remark{Remark}
We have shown that
the \cstar-envelope of
 a direct limit of full upper triangular
matrix algebras with compression embeddings
 is an AF \cstar-algebra.
If, on the other hand, $A$ is the direct limit
of full upper triangular matrix algebras with
\plstar-extendible embeddings, then $A$ can
be represented as a generating subalgebra of
a UHF algebra.  Since UHF algebras are simple,
it is immediate that the generated \cstar-algebra
is the \cstar-envelope.  More is true.  It is easy
to show that if $A$ is a  generating
subalgebra of an AF \cstar-algebra $B$
and if $A$ contains a Stratila-Voiculescu masa,
 then $B$ is
the \cstar-envelope of $A$.  Thus, if $A$ is the
direct limit of a system of digraph algebras with
\plstar-extendible embeddings, then the \cstar-envelope
of $A$ is an AF algebra.
\par
The simple proof of this last fact does not apply
to compression limits.  Examples in the next section
will show that the image of the limit algebra under
the representation $\rho$ need not contain a masa in
the generated \cstar-algebra.  (The natural diagonal
in $\rho(A)$ need not be a masa in $C^*(\rho(A))$.)
The following problem is suggested.
\endremark
\remark{Problem}
Is the \cstar-envelope of a direct limit of digraph
algebras with compression embeddings an AF \cstar-algebra?
\endremark
\remark{Remark}
In \cite{P1}, Power studies a direct limit
of digraph algebras with embeddings which are
neither \plstar-extendible nor compression
embeddings.
For systems of tri-diagonal algebras with
certain natural non-\plstar-extendible embeddings,
Power shows that the limit algebra is completely
isometrically isomorphic to a generating subalgebra
of an appropriate Bunce-Deddens algebra.
Since this latter algebra is simple, it is the
\cstar-envelope of the limit algebra for the
tri-diagonal system.  In particular, this provides
an example of a system of digraph algebras whose limit
algebra has a \cstar-envelope which is not AF.
\endremark
\remark{Remark}
Since $\rho(A)$ is weakly dense in $\algn$, $\crho$ is
an irreducible \cstar-algebra.  If it happens that
$\rho(A)$ contains a non-zero compact operator, then
Arveson's boundary theorem \cite{A2} immediately
implies that the Silov boundary for $\rho(A)$ is 0.
The boundary theorem is a deep theorem, so the argument
above could be considered more elementary.
\par
In most of the examples in the next section,
$\rho(A)$ contains non-zero compact operators.
Here is how to determine in general if $\rho(A)$
contains compact operators.  Each representation
$\rho_k$ of $T_{n_k}$ can be written as an infinite
direct sum $\sum^{\oplus} \eta_i$, where the $\eta_i$'s
are compressions.  Let $z_k$ be the number of times that
$\eta_i = id$ in the expression for $\rho_k$.  Then
$z_k \in \{\,1,2,\dots,\infty\,\}$ and the sequence
$z_k$ is decreasing.  Thus there are two possibilities:
all $z_k = \infty$, or there is  a finite integer $y$
such that $z_k = y$ for all large $k$.  In the first
case $\rho(A)$ will contain no non-zero compact operators.
In the second case, $\rho(A)$ contains finite rank operators
and $\crho$ contains all compact operators.
\par
For the first assertion, for any $a \in T_{n_k}$,
$\|a\| = \|\rho_k(a)\| = \|\rho_k(a)\|_{\text{ess}}$.
It follows from the density of $\bigcup \rho_k(T_{n_k})$
that $\|\rho(a)\| = \|\rho(a)\|_{\text{ess}}$ for all
$a \in A$.  Thus, $\rho(A)$ contains no non-zero compact
operators.  As for the second assertion, choose $k$ so
that $z_k$ is finite.  Let $v$
be the matrix unit in $T_{n_k}$ which has a 1 in the
extreme upper right hand corner and zeros elsewhere.
Then $\rho_k(v)$ has finite rank and lies in $\rho(A)$.
Since $\rho(A)$ is weakly dense in $\algn$, it has no
non-trivial reducing subspaces.  Consequently, $\crho$ is
irreducible and contains a non-zero compact operator, which
implies that it contains all compact operators.
\endremark
\head
IV. Some Examples
\endhead
\subhead A \endsubhead
Fix an integer $i$ between 1 and $n$ and consider the
system:
$$
\system n {n +k_1} {n + k_1 + k_2} {A_i}.
$$
Each embedding $\phi_n$ is given by
$$
 a \mapsto a \oplus a_{ii}I_{k_n},
$$
where $I_k$ is the $k \times k$ identity matrix.
For example,
$$
a \mapsto \bmatrix
a & & \\
& a_{ii} & \\
& & a_{ii}
\endbmatrix
\mapsto
\bmatrix
a & & & & & \\
& a_{ii} & & & & \\
& & a_{ii} & & & \\
& & & a_{ii} & & \\
& & & & a_{ii} & \\
& & & & & a_{ii}
\endbmatrix
\mapsto
$$
\par
The representation $\rho$ of the limit algebra $A_i$
will act on a Hilbert space $\H$ with othornormal
basis $\{\,e_j\,\}$ indexed by $\Bbb N \cup \{0\}$.
Let $p_k$ denote the projection on the linear
span of $\{\,e_0, \dots ,e_k\,\}$ and $\N$ the
nest consisting of the $p_k$'s.  Also, let
$\kn$ be the algebra of compact operators
which leave $\N$ invariant.
\par
It is easy to check that for each $j$, the image
of any matrix under $\rho_j$ is the sum of
a finite rank matrix and a scalar multiple of
the identity $I$.  Consequently, $\rho(A_i)
\subseteq \kn + \cbb I$.  The containment will,
in fact, be proper.
\par
For each $k$, let $f_k = p_k - p_{k-1}$.  Let
$$
B_i = \{\,s \in \kn + \cbb I \:
  \lim_{k \to \infty} f_ksf_k = f_isf_i\,\}.
$$
The image of each $\rho_k$ is clearly
contained in $B_i$, whence $\rho(A_i)
\subseteq B_i$.  To see the reverse containment,
let $b \in B_i$.  There is a compact operator $c$
in $\kn$ and a scalar $\alpha$ such that
$b = c + \alpha I$.  Observe that
$$
b_{ii} = f_ibf_i  = \lim_{k \to \infty}
  f_kbf_k = \lim_{k \to \infty}(c_{kk} + \alpha) = \alpha.
$$
(Strictly speaking, this is abuse of notation, since, for
example, $f_ibf_i$ is not a scalar but an infinite matrix
with a single non-zero entry, viz., $b_{ii}$ in the
$i^{\text{th}}$ position on the diagonal.  In cases like
this we identify the two and the meaning should
be clear.)
\par
Now, fix an integer $k > i$ such that $b_{kk}$ is
close to $\alpha$ and $p_kcp_k$ is close to $c$.
Then $a = p_kbp_k + \alpha p_k^{\perp}$ is in
$\rho(A_i)$ and is close to $b$.  Thus, $\rho(A_i)$
is dense in, and hence equal to $B_i$.
\par \noindent \line{{\bf Fact.}
The family of algebras $A_i$ are pairwise non-isomorphic.
\hfil}  \par
Indeed, the identity and the $p_k$ are the only projections
in $\kn + \cbb I$ (or in $\bh$, for that matter) which
are invariant under $\rho(A_i)$.  But $p_k \in \rho(A_i)$ if,
and only if, $k < i$.  Thus, keeping in mind that the indexing
starts with 0, there are exactly $i + 1$ invariant projections
in $A_i$.  This establishes the claim.
Note also that the \cstar-envelope for each $A_i$ is
$\K + \cbb I$.  The ``linking condition'' that the
$i^{\text{th}}$ diagonal element is equal to the limit
of the diagonal elements disappears in the passage
to the generated \cstar-algebra.
This is an example in which the diagonal of
$\rho (A)$ is not a masa in $C^*(\rho(A))$.
\subhead B \endsubhead
This time consider the system
$$
\system n {n + k_1} {n + k_1 + k_2} A
$$
with the embeddings $\phi_j$ defined by
$$
\phi_j(a) = a \oplus a_{pp}I_{k_j},
$$
where $p = n + k_1 + \dots + k_{j-1}$, i.e.,
$a_{pp}$ is the last entry on the diagonal of
the matrix $a$.
\par
The Hilbert space $\H$ and the nest $\N$ are
as before.  This time, however, $\rho(A)
= \kn + \cbb I$; the arguments are much the
same as in the previous example.  The limit
algebra $A$ now has infinitely many invariant
projections, and so is not isomorphic to any
of the $A_i$ from the previous example.
Also, the \cstar-envelope of $A$ is evidently
once again $\K + \cbb I$.
\subhead C \endsubhead
For this example, it is most convenient if the
indices for each matrix algebra $T_m$ are drawn
from the set $\{\,-m+1,\dots, -1, 0\,\}$.
The basis for $\H$ is indexed by the non-positive
integers and the $p_k$ will be the obvious infinite
rank projections.  Fix an integer $i$ such that
$-n+1 \le i \le 0$.  The embeddings in the
system
$$
\system n {n+k_1} {n+k_1+k_2} {B_i}
$$
are given by
$$
a \overset \phi_j \to \longmapsto a_{ii}I_{k_j} \oplus a.
$$
\par
Arguments analogous to the one in example A show that
$$
\rho(B_i) = \{\,s \in \kn + \cbb I \:
  \lim_{k \to -\infty}f_ksf_k = f_isf_i\,\}.
$$
Once again, $B_i \cong B_j$ if, and only if $i=j$, and
the \cstar-envelope of each $B_i$ is $\K +\cbb I$.
This is another example in which the diagonal
of $\rho(A)$ is not a masa in the \cstar-envelope.
\par
The limit algebras in examples A and C are all non-isomorphic.
All that needs to be checked is that $B_i \not\cong A_i$
for each $i$, since when $i \ne j$, $B_i$ and $A_j$ have
a different number of invariant projections.  Observe
that if $p$ is a projection in $\rho(A_i) \cap
\rho(A_i)^*$ which is invariant under $\rho(A_i)$ and
is not the identity, then for any sequence $x_j$ in
$A_i$, the set $\{\,x_jp\,\}$ is linearly dependent.
On the other hand, if $p \ne 0$ is an invariant
projection in $\rho(B_i) \cap \rho(B_i)^*$, then
there exist operators $x_1,x_2, \dotsc \in \rho(B_i)$
such that the set $\{\,x_jp\,\}$ is linearly independent.
\subhead D \endsubhead
In this example, $\H$
is a Hilbert space with basis indexed by
 $\nbb \cup \{0\}$,
 $\N = \{\,p_k\,\}$ and $\kn$ are as before, and
$\D$ denotes the algebra of diagonal matrices.
Also, $d$ will denote the conditional expectation
from $T_n$ onto $D_n$; in other words, if $a$ is
an upper triangular matrix, $d(a)$ is the diagonal
part of $a$. Note that $d$ is a direct sum
of rank one compressions.
\par
Now consider the stationary system
$$
\system 2 4 8 A
$$
in which $\phi_k \: T_{2^k} \to T_{2^{k+1}}$ is
given by $\phi_k(a) = a \oplus d(a)$.
\par
The representations $\rho_k \: T_{2^k} \to \kn + \D$
are given by $\rho_k(a) = a \oplus d(a) \oplus
d(a) \oplus \dotsb = a \oplus \sum^{\oplus}d(a)$.
The image $\rho(A)$ of the limit algebra under the
canonical representation is thus a subalgebra of
$\kn +\D$.  Once again, the image will be a proper subalgebra.
\definition{Definition}  We say that an element  $b \in \D$
is {\it periodic\/} with period $p$
if $b_{m+p} = b_m$, for all $m$.
Here, $b_m$ denotes the matrix entry $b_{mm}$.
 We say $b$ is
{\it dyadic periodic\/} is the period $p$ is a power of 2.
\define\dap{\D_{AP}(2^{\infty})}
Let $\dap$ denote the closure in the norm topology of the
dyadic periodic elements of $\D$.
\enddefinition
\define\kno{\K_{\N}^0}
\par
Let $\kno = \{\, a \in \kn \: d(a) = 0 \, \}$.
Clearly, $\rho_k(T_{2^k}) \subseteq \kno + \dap$,
for each $k$.  Thus, $\rho(A) \subseteq \kno + \dap$.
Furthermore, an argument analogous to the one used
in example A shows that we have equality:
$\rho(A) = \kno + \dap$.  The \cstar-envelope
of $A$ is therefore $\K + \dap$.
\subhead E \endsubhead
The notation is the same as in the previous example,
but the embedding is changed to
$$
a \longmapsto a \oplus dlh(a) \oplus dlh(a)
$$
where $dlh(a)$ is the diagonal last half of $a$.
(When $a$ is a $2^k \times 2^k$ matrix, $dlh(a)$
is a $2^{k-1} \times 2^{k-1}$ matrix.)
\par
 If $a$ is a matrix in
$T_{2^k}$, then the first $2^{k-1}$ diagonal
terms in $\rho_k(a)$ bear no relation to the
remaining diagonal terms  in $\rho_k(a)$.  This
difference in behaviour compared with the last example
results in a slightly different image for the limit
algebra under the canonical representation: $\rho(A)
= \kn + \dap$.  The \cstar-envelope is unchanged.
\subhead F \endsubhead
Consider the stationary system
$$
\system 248A
$$
in which $\phi_k \: T_{2^k} \to T_{2^{k+1}}$ is
given by  $\phi_k(a) = a \oplus lh(a) \oplus lh(a)$.
Here $lh(a)$ denotes compression to the last
half of $a$, resulting in a $2^{k-1} \times 2^{k-1}$
matrix when $a$ is $2^k \times 2^k$.
\redefine\S{\Cal S(2^{\infty})}
\define\sn{\Cal S_{\N}(2^{\infty})}
Let $\S$ denote the norm closure of all the dyadic
periodic matrices.  (A matrix on $\H$ is {\it periodic}
with period $p$ if it has the form
$a \oplus a \oplus \dotsb $,
where $a$ is a $p \times p$ matrix.)
The intersection of $\S$ and $\algn$ will
be denoted by $\sn$.  With this notation, the
image of the limit algebra $A$ under the canonical
representation $\rho$ is $\kn + \sn$ and the
\cstar-envelope is $\K + \S$.
\subhead G \endsubhead
We conclude by mentioning an old example, the
system
$$
\system 248A
$$
with the standard embeddings of multiplicity 2:
 $\phi_k(a) = a \oplus a$.
The canonical representation is the Smith representation
described in \cite{OP}; $\rho(A) = \sn$ and the
\cstar-envelope is the UHF algebra $\S$.

\enddocument